\documentclass[11pt]{article}
\usepackage{graphicx} 
\usepackage[utf8]{inputenc}
\usepackage[T1]{fontenc}
\usepackage[french]{babel}
\usepackage{newpxtext} 
\usepackage{microtype} 
\usepackage{amsmath,amsfonts,amssymb,amsthm} 
\usepackage{enumitem} 
\usepackage{tikz}
\usetikzlibrary{calc,patterns,angles,quotes}

\usepackage{include}

\title{Dualité étale à la Poitou-Tate pour les tores sur des variétés définies sur un corps fini}
\author{Melvyn El Kamel-{}-Meyrigne}
\date{}

\begin{document}

\maketitle

\textbf{Résumé} Soit $k$ un corps global de caractéristique $p>0$.  Notons $\Omega_k$ l'ensemble des places de $k$ et soit $S$ un sous-ensemble non vide de $\Omega_k$. On considère un schéma $\sX \ra Spec(\CO_S)$ lisse, séparé, de type fini et $\sT$ un tore défini sur $\sX$. On étudie le groupe de Tate-Shafarevich  donné par les éléments de $H^1(\sX, \sT)$ qui s'annulent dans les groupes $H^1(\sX \otimes_{\CO_S} k_v, \sT)$ pour tout $v \in S$. On établit une dualité pour $\sT$ qui généralise la dualité de Poitou-Tate classique pour les tores à des variétés définies sur un corps fini de dimension arbitraire.\\

\textbf{Abstract} Let $k$ be a global field of characteristic $p>0$. Denote $\Omega_k$ the set of places of $k$ and let $S$ be a non-empty subset of $\Omega_k$. We consider a scheme $\sX \ra Spec(\CO_S)$ smooth, separated, of finite type and $\sT$ a tori defined over $\sX$. We study the Tate-Shafarevich group given by the elements of $H^1(\sX, \sT)$ which vanish in the group  $H^1(\sX \otimes_{\CO_S} k_v, \sT)$ for all $v \in S$. We establish a Poitou-Tate duality for $\sT$ which generalise the classical Poitou-Tate duality for tori for varieties defined over a finite field of arbitrary dimension.\\

\textbf{Mots clés :} Groupes de Tate-Shafarevich, Dualité de Poitou-Tate, Tores algébriques, Cohomologie étale.


\section{Introduction}

Soit $k$ un corps global et $\Omega_k$ l'ensemble des places de $k$. On dit qu'une famille de $k-$variétés $\sV$ satisfait le principe local-global (ou principe de Hasse) si $$
\forall Z \in \sV, \prod_{v \in \Omega_k}Z(k_v) \neq \emptyset \Ra Z(k) \neq \emptyset.$$\\
Soit $G$ un groupe algébrique commutatif lisse sur $k$. Lorsque $\sV$ est l'ensemble des $G$-torseurs, la véracité du principe local-global peut être déterminée au niveau cohomologique. En effet, puisque le groupe $H^{1}(k,G)$ classifie les $G$-torseurs, le groupe de Tate-Shafarevich, qui est défini par
 $$ \Sha^{1}(k, G) \coloneqq Ker (H^{1}(k,G) \ra \prod_{v \in \Omega_k} H^{1}(k_v,G))$$
mesure l'obstruction au principe local-global pour les $G$-torseurs.\\

L'étude des groupes de Tate-Shafarevich a conduit à la découverte de nombreux résultats de dualité. L'un des plus notoires d'entre eux est le théorème de dualité de Poitou-Tate pour les groupes finis étales : soit $G$ un groupe fini étale d'ordre $n$, le cup-produit induit un accouplement parfait de groupes finis
$$
\Sha^1(k,G) \times \Sha^{2}(k, Hom(G, \mu_n)) \ra \Q /\Z.
$$

Un théorème de dualité analogue existe pour les tores : soit $T$ un tore et $\hat{T}$ son module de caractères, il existe un accouplement parfait 
$$
\Sha^1(k,T) \times \Sha^{2}(k, \hat{T}) \ra \Q /\Z.
$$
\\
qui par passage au quotient par le sous-groupe divisible maximal de $\Sha^{2}(k, \hat{T})$ induit un accouplement parfait de groupes finis
$$
\Sha^1(k,T) \times \Sha^{2}(k, \hat{T})/Div \ra \Q /\Z.
$$.

Plus récemment, d'autres résultats de dualité pour les groupes finis étales et les tores définis sur des corps de fonctions $K=k_0(C)$ d'une courbe $C$ définie sur un corps de base $k_0$ ont été obtenus : Harari et Szamuely dans \cite{HarSza16} et Colliot-Thélène et Harari dans \cite{CTHar15} pour des corps de petite dimension cohomologique puis Izquierdo dans \cite{Iz16} pour les corps $d$-locaux. 

Cependant, la stratégie employée dans ces articles ne permet pas de traiter le cas des corps de fonctions $K=k_0(X)$ lorsque $X$ est une variété de dimension arbitraire. Néanmoins, il est possible d'obtenir des résultats portant sur la cohomologie de la variété $X$ elle-même plutôt que sur celle de son corps des fonctions. À cet effet, soient $k$ un corps global, $S$ un sous-ensemble non-vide de $\Omega_k$ et notons $\CO_S$ l'anneau des $S$-entiers dans $k$. Considérons un morphisme $X \ra \sS$ séparé et de type fini et soit $\sF$ un faisceau étale sur $\sX$. Afin d'énoncer un théorème de dualité portant sur la cohomologie de $\sX$, on définit les groupes de Tate-Shafarevich par
$$ \Sha^{i}(X,\sS,\sF) \coloneqq Ker (H^{i}(\sX,\sF) \ra \prod_{v \in S} H^{i}(\sX \otimes_{\CO_S} k_v,\sF)),$$
$$ \Sha^{i}_c(X,\sS,\sF) := Ker (H^{i}_c(\sX,\sF) \ra \prod_{v \in S} H^{i}_c(X \otimes_{\CO_S} k_v,\sF)).$$

Une première avancée dans cette direction fut réalisée par Saito dans son article \cite{Sa89} où il généralise la dualité de Poitou-Tate à des variétés de dimension arbitraire. Plus précisément, si l'on considère une variété $X$ lisse, géométriquement connexe, définie sur un corps global $k$ de caractéristique $p \geq 0$ de dimension $d$ et un faisceau $\sF$ localement constant, constructible de $\Z/m \Z$-modules sur $X$ avec $m \neq p$, Saito démontre l'existence un accouplement parfait de groupes finis
$$
\Sha^i(X, Spec(k), \sF) \times \Sha^{2d+3-i}_c(X, Spec(k), \mathcal{H}om(\sF, \mu_m^{\otimes d+1})) \ra \Q / \Z.
$$

L'idée de la preuve consiste à combiner la dualité de Poitou-Tate classique pour les modules finis avec la dualité de Poincaré pour les variétés lisses. Par la suite, Geisser et Schmidt ont amélioré ce résultat dans leur article \cite{Ge17} en remplaçant $\Omega_k$ par un sous-ensemble non-vide $S$ de places contenant les places archimédiennes : ils considèrent un schéma  $\sX \ra Spec(\CO_S)$  régulier, plat, séparé et de type fini de dimension relative $d$ et établissent un accouplement parfait de groupes finis
$$
\Sha^i(\sX, \CO_S, \sF) \times \Sha^{2d+3-i}_c(\sX, Spec(CO_S), \mathcal{H}om(\sF, \Z /m \Z)(d+1)) \ra \Q / \Z.
$$

L'objectif de cet article est de démontrer une dualité pour les tores au lieu des faisceaux constructibles similaire à celle de Geisser et Schmidt en utilisant les résultats de \cite{Sa89} tout en adaptant la stratégie employée dans \cite{HarSza16}. La difficulté principale dans notre contexte réside dans le fait qu'il faille travailler avec la cohomologie étale plutôt que la cohomologie galoisienne et que l'on doive recourir à la cohomologie à support compact pour des faisceaux qui ne sont pas de torsion.\\

Introduisons maintenant le cadre de cet texte. Soit $k$ 
un corps global de caractéristique $p > 0$ et soit $S$ un ensemble non vide de places de $k$. 
Notons $\CO_{S}$ l'anneau des $S$-entiers dans $k$ et $\sS =Spec(\CO_{S})$. Soit $\sX \ra \sS$ un schéma séparé, de type fini, 
lisse de dimension relative $d$ et $\hat{\sT}$ un faisceau défini sur $\sX$ localement isomorphe à un faisceau constant libre de type fini.
Posons $\check{\sT} = \mathcal{H}om(\hat{\sT}, \Z)$ et fixons un entier $a \in \{0,1,...,d+1\}.$ Soit $\sT = \check{\sT} \otimes \Z(a)[1]$ où $\Z(a)$ désigne le $a$-ème complexe motivique. Le complexe $\tilde{\sT} \coloneqq \hat{\sT} \otimes \Z(d+1-a)[1]$ jouera le rôle du dual de $\sT.$\\

\begin{thme}
Supposons que la torsion $l$-primaire du groupe $H^{a+1}(\sX \otimes k_v,\Z (a))$ est d'exposant fini pour tout $v \in S$ et pour tout $l \neq p$ premier.
 On a un accouplement $$\Sha^{a}(\sX,\sS,\sT)\{p'\} \times \Sha^{2d+2-a}_c(\sX,\sS,\tilde{\sT})\{p'\} \ra \Q / \Z$$
entre les torsions premières à $p$ des groupes $\Sha^{a}(\sX,\sS,\sT)$ et $\Sha^{2d+2-a}_c(\sX,\sS,\tilde{\sT})$  qui, après passage au quotient par les sous-groupes divisibles maximaux, induit un accouplement parfait de groupes finis 
    $$
    \ol{\Sha^{a}(\sX,\sS,\sT)}\{p'\} \times \ol{\Sha^{2d+2-a}_c(\sX,\sS,\tilde{\sT})}\{p'\} \ra \Q / \Z\{p'\}. 
    $$
\end{thme}

Dans le cas où $a=1$, l'isomorphisme $\Z(1)[1] \simeq \gm$ montre que le faisceau $\sT$ s'identifie à un tore et on peut alors montrer que l'hypothèse sur le groupe $H^{a+1}(X_v,\Z (a) )\{l\}$ est toujours satisfaite; on obtient alors un résultat plus précis.

\begin{thme}
Supposons que le faisceau $\sT$ s'identifie à un tore. On a un accouplement parfait de groupes finis
    $$\Sha^{1}(\sX,\sS,\sT)\{p'\} \times \overline{\Sha^{2d+1}(\sX,\sS,\tilde{\sT})}\{p'\} \ra \Q / \Z\{p'\}.$$
\end{thme}


\subsection{Rappels généraux et notations}

On rassemble ici les notations utilisées dans l'article afin que l'on puisse plus facilement s'y référer :
\begin{itemize}[label={$\bullet$}]
    \item $k$ un corps global de caractéristique $p > 0$.
    \item $k_s$ une clôture séparable de $k$.
    \item $k_v$ la complétion de $k$ en une place $v$.
    \item $\CO_v$ l'anneau des entiers de $k_v$ pour une place $v$. On note $k_v^h$ le corps des fractions de l'hensélisé $\CO_v^h$ de $\CO_v$,
    \item $\sC$ une courbe lisse propre géométriquement intègre sur un corps fini $k_0$ de corps des fonctions $k$.
    \item $\sC^{(1)}$ l'ensemble des points de codimension $1$ de $\sC$. 
    \item $S$ un sous-ensemble non-vide de $\sC^{(1)}$.
     \item $\sS \coloneqq Spec (\CO_S) \coloneqq \prl_{U \subseteq \sC, \sC \setminus U \subseteq S} U$ le spectre de l'anneau des $S$-entiers.
    \item $f : \sX \ra \sS$ un schéma séparé, lisse, de type fini de dimension relative $d$.
    \item $X_v \coloneqq \sX \otimes_{\sS} k_v$ pour une place $v$ de $k$.
    \item $\hat{\sT}$ un faisceau étale sur $\sX$ localement isomorphe à un faisceau constant libre de type fini pour la topologie étale.
    \item $\check{\sT} \coloneqq \mathcal{H}om(\hat{\sT},\Z)$.
    \item  $\sT \coloneqq \check{\sT} \otimes \Z(a)[1]$.
    \item $\tilde{\sT} := \hat{\sT} \otimes \Z(d+1-a)[1]$ qui joue le rôle du dual de $\sT$.
    \item $\tilde{\sT}_t \coloneqq \hat{\sT} \otimes \Q / \Z(d+1-a)$
    \item Pour tout faisceau étale $\sF$ sur $\sX$ de $\Z/m\Z$-modules avec $m$ premier à $p$, on note $\sF(r)$ le $r$-ème twist de Tate de $\sF$.
\end{itemize}

Soit $G$ un groupe abélien. On utilisera les notations suivantes :
\begin{itemize}[label={$\bullet$}]
    \item ${}_{n}G$ : la $n$-torsion de $G$ où $n \in \N^{*}$.
    \item $G\{p\}$ : la torsion $p$-primaire de $G$.
    \item $G\{p'\} \coloneqq \bigoplus_{l \neq p}G\{l\}$ : la torsion première à $p$.
    \item $G^{(p)} \coloneqq \prl_{n \geq 0} G/p^nG$ le complété de $G$ pour la topologie $p$-adique.
     \item $\overline{G}$ : le quotient de $G$ par son sous-groupe divisible maximal.
     \item $G^D \coloneqq Hom(G, \Q / \Z)$ le dual de Pontryagyn de $G$.
\end{itemize}







\subsection{Complexes motiviques :}
Soit $X$ une variété lisse définie au-dessus d'un corps $k$. Bloch définit dans son article \cite{Blo86} les complexes motiviques $\Z(i)$ (respectivement $\Z(i)_{Zar}$) de faisceaux étales (resp. de Zariski) définis sur $X$ pour tout entier $i \geq 0$. Pour tout groupe abélien $A$, on note $A(i)$  (resp. $A(i)_{Zar}$) le complexe $A \otimes \Z(i)$  (resp. $A \otimes \Z(i)_{Zar})$. Geisser est ensuite parvenu à étendre cette construction lorsque la base est le spectre d'un anneau de Dedekind dans \cite{GE04}. Les complexes $\Z(i)$ et $\Z(i)_{Zar}$ vérifient un certain nombre de propriétés que nous utiliserons plusieurs fois par la suite; nous les listons ci-dessous :

\begin{lemme} Soit $X$ un schéma séparé lisse de type fini défini au-dessus d'un anneau de Dedekind $B$ et soient $i,j \geq 0$.
\begin{enumerate}[label=\roman*)]
    \item $\Z(0) \simeq \Z$ et $\Z(1)[1] \simeq \G_m .$ 
    \item Le complexe $\Z(i)_{Zar}$ est concentré en degrés $\leq i$. 
    \item Soit $\alpha$ la projection du site étale de Y sur le site de Zariski de $X$, on a un isomorphisme $\Q(i)_{Zar} \simeq R\alpha_{*}\Q(i)$. 
    \item Pour $n$ inversible dans $B$, on a un quasi-isomorphisme $ \Z/n\Z(i) \simeq \mu_{n}^{\otimes i}$.
    \item Si $B$ est un corps, il existe un accouplement $\Z(i) \otimes^{\mathbf{L}} \Z(j) \ra \Z(i+j)$ dans la catégorie dérivée des faisceaux étales sur $X$.
\end{enumerate}
\end{lemme}

\proof 
\begin{enumerate}[label=\roman*)]
\item Corollaire 6.4 de \cite{Blo86}.
\item Lemme 2.5 de \cite{KA12}.
\item Proposition 5.1 de \cite{GE04}.
\item Théorème 1.2, 4) de \cite{GE04}.
\item \cite{To92}.
\end{enumerate}

\subsection{Préambule sur la cohomologie à support compact}

On considère un ouvert suffisamment petit $\sU$ de $\sC$ contenant $\sS$ de sorte que l'on puisse étendre le morphisme $\sX \ra \sS$ à un morphisme $\sX_{\sU} \ra \sU$ qui est également lisse, séparé et de type fini; on a alors que $\sX \simeq \sX_{\sU} \times_{\sU} \sS$. De plus, on peut supposer que l'on puisse identifier $\hat{\sT}$ à la restriction d'un faisceau localement constant libre de type fini sur $(\sX_{\sU})_{\textit{\'et}}$ que l'on notera aussi $\hat{\sT}$ pour alléger la notation.\\

Étant donné que les faisceaux qui nous intéressent ne sont pas de torsion, il faut être précautionneux lorsque l'on manipule la cohomologie à support compact.
Notons $g_{\sU}$ le morphisme naturel $g_{\sU} : \sX_{\sU} \ra \sC$ Soit $Z$ une compactification de $g_{\sU}$ de sorte que $g_{\sU}$ se factorise par une immersion ouverte $l : \sX_{\sU} \ra Z$ et un morphisme propre $p : Z \ra \sC$; une telle compactification existe par le théorème de compactification de Nagata, cf \cite{Conrad}. Comme $\sX_{\sU}$ est un schéma régulier, on peut quitte à remplacer $Z$ par sa normalisation supposer que $Z$ est normal. Cette compactification est fixée pour le reste de l'article.\\
Pour tout ouvert $\sR \subseteq \sU$ contenant $\sS$, on note $\sX_{\sR} \coloneqq \sX_{\sU} \times_{\sU} \sR$ et $j_{\sR} : \sX_{\sR} \ra Z$ l'immersion ouverte induite par l'inclusion $ \sR \ra \sU$ composée avec le morphisme $\sX_{\sU} \ra Z$. Le diagramme commutatif suivant résume la situation :

\[\begin{tikzcd}
\sX & \sX_\sR & \sX_\sU & Z\\
	\sS & \sR & \sU & \sC \\
 {}
	\arrow[from=1-1, to=1-2]
	\arrow["v_{\sR}", from=1-2, to=1-3]
    \arrow["u_{\sR}", from=2-2, to=2-3]
    \arrow["j_{\sR}", bend left=30, from=1-2, to=1-4]
	\arrow[from=2-1, to=2-2]
	\arrow[from=2-2, to=2-3]
	\arrow["f", from=1-1, to=2-1]
	\arrow["f_\sR", from=1-2, to=2-2]
	\arrow["f_\sU", from=1-3, to=2-3]
   \arrow["i_{\sR}" below, bend right, from=2-2, to=2-4]
    \arrow["i_{\sU}", from=2-3, to=2-4]
    \arrow["p", from=1-4, to=2-4]
    \arrow["j", from=1-3, to=1-4]
\end{tikzcd}\]

Notons que pour tout ouvert $\sR$ inclus dans $\sU$ contenant $\sS$, le schéma $Z$ est également une compactification du morphisme $g_{\sR} : \sX_{\sR} \ra \sC$. On peut alors définir la cohomologie à support compact de $\sX_{\sR}$ par rapport à $\sC$ par $H^i_{c, \sC}(\sX_{\sR}, \sF) \coloneqq H^i(\sC, R{g_{\sR}}_!\sF)$ pour tout faisceau étale $\sF$ sur $\sX_{\sR}$, où $R{g_{\sR}}_! \coloneqq Rp_{*}{j_{\sR}}_!$  par définition.\\

De manière analogue, on définit pour tout ouvert $\sX_\sR$ de $\sX_\sU$ la cohomologie à support compact par rapport à $\sR$ de la façon suivante : le schéma $Z_\sR \coloneqq Z \times_\sC \sR \ra \sR$ est une compactification du morphisme $f_{\sR} \sX_\sR \ra \sR$  par changement de base et en considérant le diagramme 
\[\begin{tikzcd}
\sX_\sR & Z_\sR & Z \\
& \sR & \sC, \\
	\arrow["l_{\sR}", from=1-1, to=1-2]
	\arrow["z_{\sR}", from=1-2, to=1-3]
    \arrow["i_{\sR}", from=2-2, to=2-3]
	\arrow["f_\sR", from=1-1, to=2-2]
	\arrow["p_{\sR}", from=1-2, to=2-2]
	\arrow["p", from=1-3, to=2-3]
\end{tikzcd}\]
 on peut définir $H^i_{c, \sR}(\sX_{\sR}, \sF) \coloneqq H^i(\sR, R{f_{\sR}}_!, \sF) $ où $R{f_{\sR}}_! \coloneqq {Rp_{\sR}}_{*}{l_{\sR}}_!$ par définition.\\

 Posons également $Z_{\sS} \coloneqq Z \times_\sC \sS$. Le morphisme naturel $l_{\sS} : \sX \ra Z_{\sS}$ permet de définir le groupe $H^i_{c, \sS}(\sX, \sF) \coloneqq H^i(\sR, Rf_! \sF)$.

 Enfin, pour toute place $v$ de $k$, le schéma $Z_v \coloneqq Z \otimes_\sC k_v$ est une compactification de $X_v \coloneqq \sX \otimes_\sS k_v$ et en considérant le diagramme ci-dessous

 \[\begin{tikzcd}
 X_v & Z_v \\
 & k_v \\
	\arrow["l_{v}", from=1-1, to=1-2]
	\arrow["f_{v}", from=1-1, to=2-2]
	\arrow["p_{v}", from=1-2, to=2-2]
\end{tikzcd}\]

on pose $H^i_{c, k_v}(\sX_{v}, \sF) \coloneqq H^i(k_v, R{f_{v}}_!, \sF) $ où $R{f_{v}}_! \coloneqq R{p_{v}}_* \circ {l_{v}}_!$ par définition.\\

Même lorsque les faisceaux considérés ne sont pas de torsion, on dispose quand même de théorèmes de changement de base dans certains cas nous permettant d'obtenir les résultats suivants.

\begin{prop}\label{propcontracohomologiecompactabsolue} Covariance de la cohomologie à support compact absolue.\\  Soient $\sP \subseteq \sR$ deux ouverts de $\sU$ contenant $\sS$ et soit $\sF$ un faisceau étale abélien défini sur $\sX_{\sU}$.
Alors, l'inclusion $\sX_{\sR} \subseteq \sX_{\sU}$
induit un morphisme fonctoriel $H^i_{c, \sC}(\sX_{\sR
}, \sF) \ra H^i_{c, \sC}(\sX_{\sU}, \sF)$ de sorte que le diagramme

\begin{equation}\label{diag1}
\begin{tikzcd}
H^i_{c, \sC}(\sX_{\sP
}, \sF) & H^i_{c, \sC}(\sX_{\sR}, \sF)  \\
& H^i_{c, \sC}(\sX_{\sU
}, \sF)  \\
	\arrow[from=1-1, to=1-2]
	\arrow[from=1-2, to=2-2]
	\arrow[from=1-1, to=2-2]
\end{tikzcd}
\end{equation}

commute.
\end{prop}

\proof Notons $\iota_{\sR}$ l'immersion fermée $\sX_{\sU} \setminus \sX_{\sR} \ra \sX_{\sU}$. La suite exacte de faisceaux étales $$
0 \ra {v_{\sR}}_!{v_{\sR}}^*  \sF \ra \sF \ra {\iota_{\sR}}_! \iota_{\sR}^* \sF \ra 0
$$ donne après avoir appliqué le foncteur exact ${Rg_{\sU}}_!$ le triangle distingué 
$$
 {Rg_{\sR}}_!{v_{\sR}}^*\sF  \ra {Rg_{\sU}}_!\sF \ra {Rg_{\sU}}_!\iota_! \iota^* \sF \ra {Rg_{\sR}}_!{v_{\sR}}^*\sF[1]
$$
qui induit un morphisme  $H^i_{c, \sC}(\sX_{\sR}, \sF) \ra H^i_{c, \sC}(\sX_{\sU}, \sF)$.\\
Étant donné que le diagramme 
\[\begin{tikzcd}
R{g_{\sP}}_!  & R{g_{\sR}}_!  \\
& R{g_{\sU}}_!\\
	\arrow[from=1-1, to=1-2]
	\arrow[from=1-2, to=2-2]
	\arrow[from=1-1, to=2-2]
\end{tikzcd}\]
commute, on en déduit la commutativité du diagramme (\ref{diag1}) en passant à la cohomologie.\qed

\begin{prop}\label{propconvacohomologiecompactrelative}Contravariance de la cohomologie à support compact relative.\\
    Soient $\sP \subseteq \sR$ deux ouverts de $\sU$ contenant $\sS$ et soit $\sF$ un faisceau étale abélien défini sur $\sX_{\sU}$.
    Alors, le morphisme de restriction $H^i(\sU, {Rf_{\sU}}_! \sF) \ra H^i(\sR, {Rf_{\sR}}_! \sF)$ induit un morphisme de restriction fonctoriel $H^i_{c, \sU}(\sX_{\sU}, \sF) \ra H^i_{c, \sR}(\sX_{\sR}, \sF)$ de sorte que le diagramme

\[\begin{tikzcd}
H^i_{c, \sU}(\sX_{\sU}, \sF) & H^i_{c, \sR}(\sX_{\sR}, \sF)  \\
& H^i_{c, \sP}(\sX_{\sP}, \sF) \\
	\arrow[from=1-1, to=1-2]
	\arrow[from=1-2, to=2-2]
	\arrow[from=1-1, to=2-2]
\end{tikzcd}\]

commute.
\end{prop}

\proof Il suffit de montrer que $H^i(\sR, u_{\sR}^*{Rf_{\sR}}_! \sF) \simeq H^i_{c, \sR}(\sX_{\sR}, \sF_{\mid \sX_{\sR}})$. Considérons pour cela le diagramme commutatif suivant

 \[\begin{tikzcd}
\sX_{\sR} & \sX_{\sU} \\
Z_{\sR} & Z_\sU \\
\sR & \sU \\
	\arrow["p_{\sR}", from=2-1, to=3-1]
	\arrow["u_{\sR}", from=3-1, to=3-2]
	\arrow["p_{\sU}", from=2-2, to=3-2]
    \arrow["s_{\sR}", from=2-1, to=2-2]
    \arrow["v_{\sR}", from=1-1, to=1-2]
    \arrow["l_{\sR}", from=1-1, to=2-1]
    \arrow["l_{\sU}", from=1-2, to=2-2]
    \arrow["f_{\sR}" left, bend right=40, from=1-1, to=3-1]
    \arrow["f_{\sU}", bend left=40, from=1-2, to=3-2]
\end{tikzcd}\]

Puisque le morphisme $s_{\sR} : Z_{\sR} \ra Z_{\sU}$ est une immersion ouverte, le théorème de changement de base de \cite{De88} fournit un premier isomorphisme de foncteurs sur $Ab((Z_{\sU})_{\textit{\'et}})$
\begin{equation}\label{preuveisochgtdebase}
    u_{\sR}^*{Rp_{\sU}}_* \simeq {Rp_{\sR}}_*z_{\sR}^*.
\end{equation}

En outre, étant donné que $\sX_{\sR} = \sX_{\sU} \times_{Z_{\sU}} Z_{\sR}$, on a par \cite[\href{https://stacks.math.columbia.edu/tag/03S6}{Tag 03S6}]{stacks-project} un second isomorphisme de foncteurs sur $Ab((\sX_{\sU})_{\textit{\'et}})$
\begin{equation}\label{preuveisochgtdebase!}
    s_{\sR}^*{l_{\sU}}_! \simeq {l_{\sR}}_!v_{\sR}^*.
\end{equation}

On en déduit alors les isomorphismes suivants
\begin{align*}
H^i(\sR, u_{\sR}^*{Rf_{\sU}}_! \sF) & = H^i(\sR, u_{\sR}^*{Rp_{\sU}}_* {l_{\sU}}_! \sF) \\
& \simeq  H^i(\sR, {Rp_{\sR}}_*s_{\sR}^*{l_{\sU}}_! \sF) \hspace{3cm }\text{par l'isomorphisme (\ref{preuveisochgtdebase}).} \\
& \simeq  H^i(\sR, {Rp_{\sR}}_* {l_{\sR}}_!v_{\sR}^* \sF) \hspace{3cm }\text{par l'isomorphisme (\ref{preuveisochgtdebase!}).} \\
& = H^i_{c, \sR}(\sX_{\sR}, \sF_{\mid \sX_{\sR}}).
\end{align*}

La commutativité du diagramme s'obtient par fonctorialité des morphismes considérés. \qed

\begin{prop}\label{propcohomologiecompactcomplétion}
Soit $\sF$ un faisceau étale abélien défini sur $\sX$. Alors, le morphisme de restriction $H^i(\sS, Rf_{!}\sF) \ra H^i(k_v, {Rf_v}_{!}\sF)$ induit un morphisme fonctoriel $H^i_{c, \sS}(\sX, \sF) \ra H^i_{c, k_v}(X_v, \sF).$
\end{prop}

\proof La preuve est identique à celle de la proposition \ref{propconvacohomologiecompactrelative} en notant que le morphisme $Z_v \ra Z_{\sS}$ est plat par construction. \qed \\

\begin{prop}\label{proplienentrecohomlogiecompact}Lien entre les cohomologies à support compact.\\
Soit $\sR$ un ouvert de $\sU$ contenant $\sS$ et soit $\sF$ un faisceau étale abélien défini sur $\sX_{\sU}$.
Alors, le morphisme de restriction $H^i(\sC, {Rg_{\sR}}_! \sF) \ra H^i(\sR, {Rg_{\sR}}_! \sF)$ induit un morphisme de restriction $H^i_{c, \sC}(\sX_{\sR}, \sF) \ra H^i_{c, \sR}(\sX_{\sR}, \sF)$.
\end{prop}

\proof Étant donné que $H^i_{c, \sC}(\sX_{\sR}, \sF) = H^i(\sC, {Rg_{\sR}}_! \sF)$, il suffit de montrer que $H^i(\sR, {Rg_{\sR}}_! \sF) \simeq H^i_{c, \sR}(\sX_{\sR}, \sF)$.\\

On utilise à nouveau le théorème de changement de base de \cite{De88} pour obtenir l'isomorphisme
\begin{equation}\label{preuveisochgtdebase!2}
    i_{\sR}^*Rp_* \simeq {Rp_{\sR}}_*z_{\sR}^*.
\end{equation}

Il en résulte que
\begin{align*}
    H^i(\sR, {Rg_{\sR}}_! \sF) & = H^i(\sR, i_{\sR}^*Rp_*{j_{\sR}}_! \sF) \\
    & \simeq H^i(\sR, {Rp_{\sR}}_*z_{\sR}^*{j_{\sR}}_! \sF)  \hspace{3cm }\text{par l'isomorphisme (\ref{preuveisochgtdebase!2}).}\\
    & = H^i(\sR, {Rp_{\sR}}_!z_{\sR}^*{z_{\sR}}_!{l_{\sR}}_! \sF) \\
    & \simeq H^i(\sR, {Rp_{\sR}}_*z_{\sR}^*{l_{\sR}}_! \sF)  \hspace{3cm }\text{puisque $z_{\sR}^*{z_{\sR}}_! \simeq id$.} \\
    & = H^i_{c, \sR}(\sX_{\sR}, \sF).  
\end{align*}
\qed \\

\begin{prop}\label{proplimitecohomologie}Commutativité entre la cohomologie à support compact et les limites directes.\\
Soit $\sF$ un faisceau étale abélien défini sur $\sX_{\sU}$. Alors, on dispose d'un isomorphisme
$$
H^i_{c, \sS}(\sX, \sF) \simeq \drl_{\sS \subseteq \sR, \sR \hspace{0.1 cm}\text{ouvert de } \sU} H^i_{c, \sR}(\sX_{\sR}, \sF).
$$
\end{prop}

\proof 
Pour tout ouvert $\sR \subseteq \sU$ contenant $\sS$, on a un isomorphisme par changement de base  (\cite[\href{https://stacks.math.columbia.edu/tag/03S6}{Tag 03S6}]{stacks-project}) $${u_{\sR}}^*{l_{\sU}}_! \sF \simeq {l_{\sR}}_!{v_{\sR}}^* \sF.$$ 
De plus, on a un isomorphisme 
$$({l_{\sU}}_! \sF)_{\mid \sS}  \simeq l_!(\sF)_{\mid \sS}.$$
également par changement de base.\\

On peut alors utiliser le corollaire 5.8, exposé $\textit{VII},$ de \cite{SGA4} pour prouver que
$$
H^i(Z_\sS, l_!  \sF_{\mid \sS})) \simeq \drl_{\sS \subseteq \sR, \sR \hspace{0.1 cm}\text{ouvert}} H^i(Z_{\sR}, {l_{\sR}}_! (\sF_{\mid \sR})),
$$
c'est ce que l'on voulait. \qed \\

\begin{core}\label{cormorphismerestrictionlimite}
    Soit $\sR$ un ouvert de $\sU$ contenant $\sS$. Alors, il existe un morphisme de restriction fonctoriel $H^i_{c, \sR}(\sX_{\sR}, \sF) \ra H^i_{c, \sS}(\sX, \sF)$.
\end{core}

\proof C'est une conséquence directe des propositions \ref{propconvacohomologiecompactrelative} et \ref{proplimitecohomologie} en passant à la limite.\qed \\

\section{Une dualité à la Artin Verdier}\label{sectiondualitéglobale}

L'accouplement $\Z(i) \otimes^{\mathbf{L}} \Z(j) \ra \Z(i+j)$ dans la catégorie dérivée des faisceaux étales sur $\sX_{\sU}$ induit l'accouplement 
\begin{equation*}
(\check{\sT} \otimes^{\mathbf{L}} \Z(a)[1]) \otimes^{\mathbf{L}} (\hat{\sT} \otimes^{\mathbf{L}} \Z(d+1-a)[1]) \ra \Z(d+1)[2].
\end{equation*}

Passer à la cohomologie fournit pour tout ouvert $\sR$ de $\sU$ les accouplements
\begin{equation}\label{accouplementcohomologieglobale}
    H^{i}(\sX_{\sR},\sT) \times H^{2d+2-i}_{c, \sC}(\sX_{\sR},\tilde{\sT})   \ra H ^{2d+4}_{c, \sC}(\sX_{\sR},\Z (d+1)).
\end{equation}

\begin{prop}\label{propQ/Zglobale}
 Pour tout ouvert $\sR$ de $\sU$ contenant $\sS$, on a un isomorphisme $H^{2d+4}_{c, \sC}(\sX_{\sR},\Z (d+1))\{p'\} \simeq \Q / \Z\{p'\}$
\end{prop}

\proof
La suite exacte de complexes $$ 0 \ra \Z(d+1) \ra \Q(d+1) \ra \Q / \Z(d+1) \ra 0$$ induit , après avoir appliqué le foncteur exact $R{g_{\sR}}_!$,
une suite exacte cohomologique $$H^{2d+3}_{c, \sC}(\sX_{\sR},\Q (d+1)) \ra H^{2d+3}_{c, \sC}(\sX_{\sR},\Q / \Z (d+1)) \ra H^{2d+4}_{c, \sC}(\sX_{\sR},\Z (d+1)) \ra H^{2d+4}_{c, \sC}(\sX_{\sR},\Q (d+1)). $$ Notons $\alpha$ la projection du petit site étale de $X$ sur son petit site de Zariski. D'après (\cite{KA12}, Theorem 2.6, c),  le morphisme naturel $R{g_{\sR}}_!\Q(d+1)_{Zar} \ra R\alpha_{*}R{g_{\sR}}_!\Q(d+1)$ est un isomorphisme dans la catégorie dérivée des faisceaux de Zariski. Ainsi, pour tout $i>0$, $$
H^{i}_{c, \sC}(\sX_{\sR},\Q (d+1)) \simeq H^{i}(Z,{j_{\sR}}_!\Q (d+1)_{Zar}) \simeq H^{i}_{Zar}(Z,{j_{\sR}}_!\Q (d+1)_{Zar}) 
$$


et ce dernier groupe est nul dès que $i>2d+2$ puisque le complexe $j_!\Q(d+1)_{Zar}$ est concentré en degré au plus $d+1$ et que la variété $Z$ est de même dimension que $\sX_{\sR}$, c'est-à-dire de dimension $d+1$. Il en découle un isomorphisme $$ H^{2d+4}_{c, \sC}(\sX_{\sR},\Z (d+1)) \simeq H^{2d+3}_{c, \sC}(\sX_{\sR},\Q / \Z (d+1)) \simeq H^{2d+3}(Z,{j_{\sR}}_!\Q / \Z (d+1)). $$

Notons $(k_0)_s$ une clôture séparable de $k_0$ et posons $Z_s \coloneqq Z \times_{k_0} ({k_0})_s$,  la suite spectrale de Hochschild-Serre $H^{i}(k_0,H^{j}(Z_s,{j_{\sR}}_!\mu_{m}^{\otimes d+1})) \Rightarrow H^{i+j}(Z, {j_{\sR}}_!\mu_{m}^{\otimes d+1})$ dégénère dès que $i>1$ ou que $j>2d+2$ puisque les dimensions cohomologiques pour la topologie étale de $k_0$ et de $Z_s$ sont respectivement $1$ et $2d+2$. On en déduit les isomorphismes 

\begin{align*}
H^{2d+3}(Z,{j_{\sR}}_!\Q / \Z (d+1))\{p'\} & \simeq \underset{p \nmid m}\drl \; H^{1}(k_0,H^{2d+2} (Z_s,{j_{\sR}}_!\mu_{m}^{\otimes d+1}))\\
& \simeq \underset{p \nmid m}\drl \; H^{1}(k_0,\Z/m\Z) \hspace{1cm}\text{par le lemme 11.3 de \cite{Mi80}} \\
& \simeq \Q / \Z\{p'\},
\end{align*}
ce qui permet de conclure. \qed \\

\rmke L'isomorphisme $H^{2d+2} (Z_s,{j_{\sR}}_!\mu_{m}^{\otimes d+1}) \simeq \Z/m\Z$ qui provient de la dualité de Poincaré est vrai même lorsque $Z$ n'est pas lisse, cela nous permet de ne pas avoir à faire l'hypothèse supplémentaire que $\sX$ possède une compactification lisse qui n'est pas garantie en caractéristique $p>0$.\\

Nous aurons d'abord besoin d'une dualité à la Artin-Verdier pour les faisceaux constructibles localement constants afin de pouvoir établir une dualité du même type pour les tores; nous en donnons l'énoncé ci-dessous.

\begin{prop}\label{propdualitéglobalefinis} \cite{Mi06}, \textit{II}, corollaire 7.7.\\
    Soit $\sR$ un ouvert de $\sU$ contenant $\sS$ et soit $\sF$ un complexe de faisceaux constructibles localement constants de $\Z /m\Z$-modules défini sur $(\sX_{\sR})_{\textit{\'et}}$ avec $m$ premier à $p$. Posons $\sF^{\vee} \coloneqq  \mathcal{H}om(\sF, \Q / \Z)$. Alors, le cup-produit induit un accouplement parfait de groupes finis
    $$
    H^{i}(\sX_{\sR},\sF) \times H^{2d+3-i}_{c, \sC}(\sX_{\sR},\sF^{\vee}(d+1)) \ra \Z /m\Z
    $$
    pour tout $i \in \Z$.
\end{prop}

On peut maintenant établir une dualité partielle pour les faisceaux $\sT$ et $\tilde{\sT}$.

\begin{prop}\label{propdualitéglobale}
   Pour tout ouvert $\sR$ de $\sU$ contenant $\sS$, l'accouplement \ref{accouplementcohomologieglobale} induit pour tout nombre premier $l \neq p$ des accouplements parfaits de groupes finis 
    $$ (H^{i}(U,\sT)\{l\})^{(l)} \times H^{2d+2-i}_{c, \sC}(U, \tilde{\sT})^{(l)}\{l\}   \ra \Q_l / \Z_l $$
\end{prop}

\proof On emploie le même raisonnement que dans \cite{HarSza08}, théorème 1.3, où l'on déduit le résultat voulu de la dualité déjà connue lorsque les coefficients sont finis. On redonne l'argument ci-dessous.\\

Le triangle distingué $$ \check{\sT}\otimes \Z / l^n\Z(a) \ra \sT \ra \sT \ra \check{\sT} \otimes \Z / l^n\Z(a)[1] $$
fournit une suite exacte  de cohomologie $$
0 \ra H^{i-1}(U,\sT)/l^n \ra H^{i}(U,\check{\sT} \otimes \Z / l^n\Z(a)) \ra {}_{l^n}H^{i}(U,\sT) \ra 0. $$

Passer à la limite inductive donne la suite exacte
$$
0 \ra H^{i-1}(U,\sT) \otimes \Q_l / \Z_l  \ra \lim_{\underset{n}{\ra}} H^{i}(U,\check{\sT} \otimes \Z / l^n\Z(a)) \ra H^{i}(U,\sT)\{l\} \ra 0. $$

La divisibilité du groupe $H^{i-1}(U,\sT) \otimes \Q_l / \Z_l$ entraîne pour chaque $m>0$ un isomorphisme $$(\drl_{n} H^{i}(U,\check{\sT} \otimes \Z / l^n\Z(a))) / l^m \simeq (H^{i}(U,\sT)\{l\}) / l^m .$$

On peut ensuite passer à la limite projective pour obtenir un premier isomorphisme
\begin{equation}\label{isopreuvedualitéArtin1}
(\drl_{n} H^{i}(U,\check{\sT} \otimes \Z / l^n\Z(a)))^{(l)} \simeq (H^{i}(U,\sT)\{l\})^{(l)}.
\end{equation}

Par ailleurs, le triangle distingué
$$  \hat{\sT} \otimes \Z / l^n\Z(d+1-a) \ra \tilde{\sT} \ra \tilde{\sT} \ra \hat{\sT} \otimes \Z / l^n\Z(d+1-a)[1]$$ 
permet d'obtenir la suite exacte de cohomologie à support compact
$$ 0 \ra H^{2d+2-i}_{c, \sC}(U,\tilde{\sT} )/l^n \ra H^{2d+3-i}_{c, \sC}(U, \hat{\sT} \otimes \Z / l^n\Z(a)) \ra {}_{l^n}H^{2d+3-i}_{c, \sC}(U,\tilde{\sT}) \ra 0. $$

L'application $H^{2d+2-i}_{c, \sC}(U,\tilde{\sT})/l^{n+1} \ra H^{2d+2-i}_{c, \sC}(U,\tilde{\sT})/l^n$ étant surjective pour tout entier $n$, passer à la limite projective préserve l'exactitude (\cite{stacks-project},10.86, Mittag-Leffler systems); on obtient alors la suite exacte $$0 \ra H^{2d+2-i}_{c, \sC}(U,\tilde{\sT})^{(l)} \ra \prl_{n} H^{2d+3-i}_{c, \sC}(U,\hat{\sT} \otimes \Z / l^n\Z) \ra \prl_{n} ({}_{l^n}H^{2d+3-i}_{c, \sC}(U,\tilde{\sT})) \ra 0.$$

Le groupe de droite étant sans torsion, il en découle un second isomorphisme 
\begin{equation}\label{isopreuvedualitéArtin2}
    H^{2d+2-i}_{c, \sC}(U,\tilde{\sT})^{(l)}\{l\} \simeq \prl_{n}(H^{2d+3-i}_{c, \sC}(U, \hat{\sT} \otimes \Z / l^n\Z(d+1-a))\{l\}.
\end{equation}

On est donc ramenés à établir un accouplement parfait entre les groupes $\displaystyle (\drl_{n}H^{i}(U,\check{\sT} \otimes \Z / l^n\Z(a)))^{(l)}$ et $\displaystyle \prl_{n}(H^{2d+3-i}_{c, \sC}(U, \hat{\sT} \otimes \Z / l^n\Z(d+1-a))\{l\}$ d'après les isomorphismes (\ref{isopreuvedualitéArtin1}) et (\ref{isopreuvedualitéArtin2}). La proposition V.I.I. de \cite{Mi06} permet de représenter le faisceau localement constant à tiges finies $\check{\sT} \otimes \Z / l^n\Z(a)$ par un groupe fini étale. Ainsi, 
$$ 
\mathcal{H}om(\check{\sT} \otimes \mu_{l^n}^{\otimes a}, \mu_{l^n}^{\otimes d+1}) \simeq \mathcal{H}om(\check{\sT}, \mu_{l^n}^{\otimes d+1-a})  \simeq \hat{\sT} \otimes \Z / l^n\Z(d+1-a).$$

Pour tout $n>0$, les groupes $H^{i}(U,\check{\sT} \otimes \Z / l^n\Z(a))$ et $H^{2d+3-i}_{c}(U, \hat{\sT} \otimes \Z / l^n\Z(d+1-a))$ sont duaux d'après la proposition \ref{propdualitéglobalefinis}; cela nous permet de conclure en passant à la limite. \qed \\

\rmke On peut également démontrer qu'il existe un accouplement
$$ (H^{i}(U,\tilde{\sT})\{l\})^{(l)} \times H^{2d+2-i}_{c, \sC}(U,\sT)^{(l)}\{l\}   \ra \Q_l / \Z_l $$
parfait de groupes finis.

\section{Dualité locale}
On établit ici une dualité locale en suivant la même approche que dans la section \ref{sectiondualitéglobale}. Considérons une place $v \in S$ et notons $X_v \coloneqq \sX \otimes_{\sS} k_v$. Le morphisme $\sT \otimes \tilde{\sT} \ra \Z(d+1)[2]$ induit à nouveau un accouplement 
\begin{equation}
    H^{i}(X_v ,\sT) \times H^{2d+1-i}_{c,k_v}(X_v ,\tilde{\sT})   \ra
H^{2d+3}_{c, k_v}(X_v ,\Z(d+1)).
\end{equation}

\begin{prop}
    Soit $v \in S$,  alors $H^{2d+3}_{c, k_v}(X_v,\Z (d+1))\{p'\} \simeq \Q / \Z\{p'\}$.
\end{prop}


\proof La preuve est similaire à celle de la proposition \ref{propQ/Zglobale} : pour tout $v \in S$, le groupe $H^{i}(Z_v,{l_{v}}_!\Q (d+1))$ est nul dès que $i > d+1$ puisque le schéma $Z_v$ est de dimension $d$. On en déduit un isomorphisme $$ H^{2d+2}(Z_v,{l_{v}}_!\Q / \Z (d+1)) \simeq H^{2d+3}(Z_v, {l_{v}}_!\Z (d+1)). $$ Posons $\overline{Z_v} : = Z_v \times_{k_v} (k_v)^s;$ la suite spectrale de Hochschild-Serre $H^{i}(k_v,(H^{j}(\overline{Z_v},{l_{v}}_!\mu_{m}^{\otimes d+1})) \Rightarrow H^{i+j}(Z_v, {l_{v}}_!\mu_{m}^{\otimes d+1})$ dégénère dès que $i>2$ ou que $j>2d$. On en déduit des isomorphismes $$ H^{2d+2}_{c, k_v}(Z_v,\Q / \Z (d+1))\{p'\} \simeq \underset{p \nmid m}\drl H^{2}(k_v,H^{2d}_{c, k_v} (\overline{Z_v},\mu_{m}^{\otimes d+1})) \simeq \underset{p \nmid m}\drl \; H^{2}(k_v,\Z/m\Z(1)) \simeq \Q / \Z\{p'\} $$ ce qui permet de conclure. \qed \\

Comme dans le cas global, nous aurons besoin d'une dualité pour les faisceaux constructibles localement constants.

\begin{prop}\label{propdualitélocalefinis}
    Soit $v \in S$ et soit $\sF$ un faisceau constructible localement constant de $\Z /m\Z$-modules défini sur $(\sX)_{\textit{\'et}}$ avec $m$ premier à $p$. Posons $\sF^{\vee} \coloneqq \mathcal{H}om(\sF, \Q / \Z)$. Alors, il existe un accouplement parfait de groupes finis
    
\begin{equation}\label{accouplementlocal}
    H^{i}(X_{k_v},\sF) \times H^{2d+2-i}_{c, k_v}(X_{k_v},\sF^{\vee}(d+1)) \ra \Z /m\Z
\end{equation}
pour tout $i \in \Z$.
\end{prop}

\proof D'après (\cite{Sa89}, théorème 2.1, 3), on a un accouplement parfait de groupes finis
$$
H^{i}(k_v, {Rf_v}_!\sF) \times H^{2-i}_{c, k_v}(k_v,R\mathcal{H}om_{k_v, \Z /m\Z }({Rf_v}_!\sF, \Z /m\Z(1))) \ra \Z /m\Z
$$

La dualité de Poincaré (\cite{SGA4}, \textit{XVIII}, théorème 3.2.5) donne un isomorphisme 
$$
R\mathcal{H}om_{\sR, \Z /m\Z }({Rf_v}_!\sF, \Z /m\Z(1)) \simeq{Rf_v}_*R\mathcal{H}om_{k_v, \Z /m\Z }(\sF, \Z /m\Z(d+1))[2d]
$$

Comme $\Z /m \Z$ est un $\Z /m \Z$-module injectif, on a que $R\mathcal{H}om_{\sR, \Z /m\Z }(\sF, \Z /m\Z(d+1)) \simeq \sF^{\vee}(d+1)$ et on en déduit l'accouplement voulu. 
\qed

\begin{prop}\label{propdualitélocale}
    Pour tout nombre premier $l \neq p$, l'accouplement (\ref{accouplementlocal}) induit un accouplement parfait de groupes finis $$ (H^{i}(X_v,\sT)\{l\})^{(l)} \times H^{2d+1-i}_{c, k_v}(X_v,\tilde{\sT})^{(l)}\{l\}   \ra \Q_l / \Z_l$$
\end{prop}

\proof La preuve est identique à celle de la proposition \ref{propdualitéglobale} mais se base sur la dualité de la \ref{propdualitélocalefinis}. \qed \\



\section{Dualité pour les groupes de Tate-Shafarevich}

On souhaite maintenant combiner ces deux dualités partielles afin d'en déduire une nouvelle au niveau des groupes de Tate-Shafarevich
$$ \Sha^{i}(\sX,\sS,\sT) \coloneqq Ker (H^{i}(\sX,\sT) \ra \prod_{v \in S} H^{i}(\sX \otimes_{\sS} k_v,\sT))$$
et
$$ \Sha^{i}_c(\sX,\sS,\sT) := Ker (H^{i}_c(\sX,\sT) \ra \prod_{v \in S} H^{i}_c(\sX \otimes_{\sS} k_v,\sT)).$$
On démontre à cet effet quelques résultats sur les groupes cohomologiques impliqués.

\begin{lemme}\label{lemnulQ(j)}
Soit $f : Y \ra B$ un schéma normal de dimension $y$ défini sur un schéma de Dedekind $B$ et soit $\sF$ un faisceau sur $Y$ localement isomorphe à un faisceau constant libre de rang fini.\\
Soit $i \geq j + y + 1$. Alors les groupes $H^i(Y, \sF \otimes \Q(j))$ et $H^i_{c,B}(Y, \sF \otimes \Q(j))$ sont triviaux pour n'importe quelle compactification de $f : Y \ra B$.
\end{lemme}

\proof La nullité des groupes $H^i(Y, \sF \otimes \Q(j))$ se démontre de la même manière que la preuve du lemme 3.1 de \cite{Iz16} : on considère un morphisme fini étale $Y' \ra Y$ déployant $\sF$ et on note $n$ le degré de ce morphisme. D'après (\cite{KA12}, Theorem 2.6, c), on a un isomorphisme $H^i(Y', \Q(j)) \simeq H^i(Y',\Q(j)_{\textit{Zar}})_{\textit{Zar}}$ et comme le complexe $\Q(j)_{\textit{Zar}}$ est concentré en degré $\leq j$, le groupe $H^i(Y',\Q(j)_{\textit{Zar}})_{\textit{Zar}}$ est trivial dès que $i \geq j + y + 1$. Ainsi, le groupe $H^i(Y, \sF \otimes \Q(j))$ est d'exposant fini par restriction-corestriction et puisqu'il est également uniquement divisible, il est donc trivial.\\
Le même raisonnement s'applique pour la cohomologie à support compact :
soit $g :Y \ra Z$ une compactification de $Y$ et $g' : Y' \ra Z'$ la compactification induite par $g$ de $Y'$, alors $H^i_{c,B}(Y',\Q(j)) \simeq H^i(Z',g_!(\Q(j)_{\textit{Zar}}))_{\textit{Zar}} = 0$ ce qui implique que le groupe $H^i_{c,B}(Y, \sF \otimes \Q(j)) = H^i(Z,g_!(\sF \otimes \Q(j)))$ est d'exposant fini et est donc nul.\qed \\




Posons $\tilde{\sT}_t \coloneqq \hat{\sT} \otimes \Q / \Z(d+1-a)$ et considérons la suite exacte  $$0 \ra \Z(d+1-a) \ra \Q(d+1-a) \ra \Q / \Z (d+1-a) \ra 0.$$ Comme $\tilde{\sT}_t$ est un faisceau de torsion, il est plus simple à manipuler que le faisceau $\tilde{\sT}$ car on peut utiliser les théorèmes de changement base usuels sans restriction.\\
Par platitude du faisceau $\hat{\sT}$, cette suite induit un triangle distingué

\begin{equation}\label{triangledisctinguétores}
    \tilde{\sT}_t \ra \tilde{\sT} \ra \hat{\sT} \otimes \Q(d+1-a)[1] \ra \tilde{\sT}_t[1].
\end{equation}

\begin{lemme}\label{lemtypecofini}
Soit $\sR$ un ouvert de $\sU$ contenant $\sS$ et soit $v \in S$.
    \begin{enumerate}
        \item Pour tout $i \geq 0$, les groupes $H^{i}(\sX_\sR,\tilde{\sT}_t)\{p'\}, H^{i}_{c, \sC}(\sX_\sR,\tilde{\sT}_t)\{p'\}$ et $H^{i}_{c, k_v}(X_v,\tilde{\sT_t})\{p'\}$ sont de torsion et de type cofini.
        \item Les groupes $H^{2d+2-a}_{c, \sC}(\sX_\sR,\tilde{\sT})\{p'\}$ et $H^{2d+1-a}_{c, k_v}(X_v,\tilde{\sT})\{p'\}$ sont de torsion et de type cofini. 
     \end{enumerate}
\end{lemme}

\proof
\begin{enumerate}
    \item D'une part, étant donné que le faisceau $\tilde{\sT}_t$ est de torsion, il en va de même des groupes $H^{i}(\sX_\sR,\tilde{\sT}_t)$ et  $H^{i}_{c, \sC}(\sX_\sR,\tilde{\sT}_t)$. D'autre part, la suite de Kummer donne par platitude du faisceau $\hat{\sT}$ la suite exacte
 \begin{equation}\label{suiteKummer}
0 \ra {}_{n}\tilde{\sT}_t \ra \tilde{\sT}_t \overset{\times n}{\ra} \tilde{\sT}_t \ra 0
 \end{equation}   
    qui induit une surjection $$H^{i}(\sX_\sR,\hat{\sT} \otimes \Z / n\Z(d+1-a)) \ra {}_{n}H^{i}(\sX_\sR,\tilde{\sT}_t)$$ pour tout entier $n$ premier à $p$. Le premier groupe étant fini, le groupe $H^{i}(\sX_\sR,\sT \otimes \Q / \Z(d+1-a))\{p'\}$ est de type cofini. Le même raisonnement s'applique au groupe $H^{i}_{c, \sC}(\sX_\sR,\tilde{\sT}_t)$ : on applique le foncteur exact ${j_{\sR}}_!$ à la suite (\ref{suiteKummer}) pour obtenir un morphisme surjectif
    $$H^{i}(Z_\sR,{j_{\sR}}_!(\hat{\sT} \otimes \Z / n\Z(d+1-a))) \ra {}_{n}H^{i}(Z_\sR,{j_{\sR}}_!(\tilde{\sT}_t)).$$
    La même preuve s'applique pour le groupe $H^{i}_{c, k_v}(X_v,\tilde{\sT_t})\{p'\}$ .\\
    
    \item 
Appliquer le foncteur exact ${Rg_{\sR}}_!$ au triangle distingué (\ref{triangledisctinguétores})
donne en passant à la cohomologie la suite exacte
    $$  H^{2d+2-a}_{c, \sC}(\sX_\sR,\tilde{\sT}_t) \ra H^{2d+2-a}_{c, \sC}(\sX_\sR,\tilde{\sT}) \ra  H^{2d+3-a}_{c, \sC}(\sX_\sR,\hat{\sT} \otimes \Q(d+1-a)) .$$

Le groupe $H^{2d+3-a}(\sX_\sR,\hat{\sT} \otimes \Q(d+1-a))$ est nul par le lemme \ref{lemnulQ(j)}. La surjection qui en découle montre que le groupe $H^{2d+2-a}(\sX_\sR,\tilde{\sT})\{p'\}$ est de torsion et de type cofini. Le même argument s'applique au groupe $H^{2d+1-a}(X_v,\tilde{\sT})\{p'\}$ en utilisant la nullité du groupe  $H^{2d+2-a}(X_v,\hat{\sT} \otimes \Q(d+1-a))$. \qed \\
\end{enumerate}

Soit $\sR$ un ouvert de $\sU$ contenant $\sS$. D'après le corollaire \ref{cormorphismerestrictionlimite} et la proposition \ref{propcohomologiecompactcomplétion}, il existe des morphismes de restriction $$H^{i}_{c,\sR}(\sX_{\sR}, \tilde{\sT}) \ra H^{i}_{c, \sS}(\sX, \tilde{\sT})$$ et 
$$ H^{i}_{c, \sS}(\sX, \tilde{\sT}) \ra H^{i}_{c, k_v}(X_v, \tilde{\sT})$$
pour tout $v \in S$. La composition de ces morphismes permet de définir un morphisme $$
H^{2d+2-a}_{c, \sR}(\sX_{\sR}, \tilde{\sT}) \ra \bigoplus_{v \in R} H^{2d+2-a}_{c, k_v}(X_v, \tilde{\sT}).
    $$

\begin{lemme}\label{lemsuitelocalisation}
    Soit $\sR$ un ouvert de $\sU$ contenant $\sS$. Alors, la suite suivante
    $$
 H^{2d+2-a}_{c,\sC}(\sX_{\sR}, \tilde{\sT}) \ra H^{2d+2-a}_{c, \sR}(\sX_{\sR}, \tilde{\sT}) \ra \bigoplus_{v \in R} H^{2d+2-a}_{c, k_v}(X_v, \tilde{\sT})
    $$
    est exacte.
\end{lemme}

\proof D'après la proposition 3.1) 1, de \cite{HarSza16}, on dispose pour tout faisceau $\sF$ sur $\sR$ d'une suite exacte
$$
... \ra H^i_{c}(\sR,\sF) \ra H^i(\sR,\sF) \ra \bigoplus_{v \in R} H^i(k_v^h, \sF) \ra ...
$$
qui commence au degré $1$. Appliquer cette suite pour le faisceau $\sF \coloneqq {Rf_{\sR}}_{!}\tilde{\sT}_t$ fournit la suite exacte
$$
... \ra H^i(\sC, {i_{\sR}}_{!}{Rf_{\sR}}_{!}\tilde{\sT}_t) \ra H^i_{c, \sR}(\sX_{\sR},\tilde{\sT}_t) \ra \bigoplus_{v \in R} H^i(\sX \otimes k_v^h, {Rf_{\sR}}_{!}\tilde{\sT}_t) \ra ...
$$

Par (\cite{SGA4}, 5.1.3.1) et le lemme 5.1.6)  le morphisme naturel ${i_{\sR}}_{!}{Rp_{\sR}}_{*} \ra {Rp_{\sR}}_{*}{z_{\sR}}_{!}$ est un isomorphisme. Ainsi,
$$
H^i(\sC, {i_{\sR}}_{!}{Rf_{\sR}}_{!}\tilde{\sT}_t) =  H^i(\sC, {i_{\sR}}_{!}{Rp_{\sR}}_{!}{l_{\sR}}_{!}\tilde{\sT}_t) \simeq H^i(\sR, Rp_{*}{z_{\sR}}_{!}{l_{\sR}}_{!}\tilde{\sT}_t) =  H^{i}_{c,\sC}(\sX_{\sR}, \tilde{\sT}_t).
$$

Par la proposition 10.1 de \cite{Ge17}, la flèche $H^i(\sX \otimes k_v^h, {Rf_{\sR}}_{!}\tilde{\sT}_t) \ra H^i(X_v, {Rf_{\sR}}_{!}\tilde{\sT}_t)$ est un isomorphisme et on peut alors montrer que $H^i(X_v, {Rf_{\sR}}_{!}\tilde{\sT}_t) \simeq H^{i}_{c, k_v}(X_v, \tilde{\sT}_t)$ de la même manière à l'aide du théorème de changement de base.\\

D'autre part, le lemme \ref{lemnulQ(j)} montre que les groupes
\begin{itemize}
    \item $H^{2d+3-a}_{c, \sC}(\sX_{\sR}, \hat{\sT} \otimes \Q(d+1-a)),$
    \item $H^{2d+3-a}_{c, \sR}(\sX_{\sR}, \hat{\sT} \otimes \Q(d+1-a)),$
    \item $H^{2d+2-a}_{c, k_v}(X_v,  \hat{\sT} \otimes \Q(d+1-a)),$
    \item $H^{2d+3-a}_{c, k_v}(X_v, \hat{\sT} \otimes \Q(d+1-a)),$
\end{itemize}
sont triviaux. Il en découle le diagramme commutatif suivant dont la première ligne et les colonnes sont exactes 

\[\begin{tikzcd}
	&& 0 \\
	H^{2d+2-a}_{c,\sC}(\sX_{\sR}, \tilde{\sT}) &  H^{2d+2-a}_{c, \sR}(\sX_{\sR},\tilde{\sT}_t) & \bigoplus_{v \in R} H^{2d+2-a}_{c, k_v}(X_v, \tilde{\sT}_t) & {} \\
	H^{2d+2-a}_{c,\sC}(\sX_{\sR}, \tilde{\sT}) & H^{2d+2-a}_{c, \sR}(\sX_{\sR}, \tilde{\sT}) &\bigoplus_{v \in R} H^{2d+2-a}_{c, k_v}(X_v, \tilde{\sT}) \\
	0 & 0 & 0.
	\arrow[from=2-2, to=2-3]
	\arrow[from=2-1, to=2-2]
	\arrow[from=2-1, to=3-1]
	\arrow[from=2-2, to=3-2]
	\arrow[from=3-2, to=4-2]
	\arrow[from=2-3, to=3-3]
	\arrow[from=3-2, to=3-3]
    \arrow[from=3-1, to=4-1]
	\arrow[from=3-1, to=3-2]
	\arrow[from=3-3, to=4-3]
	\arrow[from=1-3, to=2-3]
\end{tikzcd}\]L'exactitude de la seconde suite horizontale se déduit alors de la première via une chasse au diagramme. \qed \\


Pour tout ouvert $\sR$ de $\sU$ contenant $\sS$, on définit $\mathcal{D}^{2d+2-a}_{c}(\sX_\sR,\tilde{\sT})=Im (H^{2d+2-a}_{c, \sC}(\sX_\sR,\tilde{\sT}) \ra H^{2d+2-a}_{c, \sS}(\sX,\tilde{\sT})) .$ Ce dernier morphisme est défini par la composition des morphismes 
$$
H^{2d+2-a}_{c,\sC}(\sX_\sR, \tilde{\sT}) \ra H^{2d+2-a}_{c,\sR}(\sX_\sR,\tilde{\sT}) \ra H^{2d+2-a}_{c,\sS}(\sX, \tilde{\sT})
$$
qui proviennent des propositions \ref{proplienentrecohomlogiecompact} et \ref{propconvacohomologiecompactrelative}.
\begin{prop}\label{propcomutativitédiagrammeD}
    Soient $\sP \subseteq \sR$ deux ouverts de $\sU$ contenant $\sS$. Alors, le diagramme 
\[\begin{tikzcd}
H^i_{c, \sC}(\sX_{\sP}, \tilde{\sT}) & H^i_{c, \sC}(\sX_{\sR}, \tilde{\sT})  \\
H^i_{c, \sP}(\sX_{\sP}, \tilde{\sT}) & H^i_{c, \sR}(\sX_{\sR}, \tilde{\sT})  \\
	\arrow[from=1-1, to=1-2]
	\arrow[from=1-2, to=2-2]
	\arrow[from=1-1, to=2-1]
    \arrow[from=2-2, to=2-1]
\end{tikzcd}\]
commute pour tout $i \geq 0.$\\

En particulier, on a l'inclusion $\mathcal{D}^{2d+2-a}_{c}(\sX_\sP,\tilde{\sT}) \subseteq \mathcal{D}^{2d+2-a}_{c}(\sX_\sR,\tilde{\sT}).$

\end{prop}

\proof
Précisons que la flèche horizontale du haut provient de la proposition \ref{propcontracohomologiecompactabsolue}.\\
D'après la proposition 3.1, 3) de \cite{HarSza16}, le diagramme 
\[\begin{tikzcd}
H^i_{c}(\sP, {Rf_{\sP}}_!\sF) & H^i_{c}(\sR, {Rf_{\sR}}_!\sF)  \\
H^i(\sP, {Rf_{\sP}}_!\sF) & H^i(\sR, {Rf_{\sR}}_!\sF)  \\
	\arrow[from=1-1, to=1-2]
	\arrow[from=1-2, to=2-2]
	\arrow[from=1-1, to=2-1]
    \arrow[from=2-2, to=2-1]
\end{tikzcd}\]
commute pour tout faisceau étale $\sF$ sur $\sX_{\sR}$. En prenant $\sF \coloneqq  \tilde{\sT}_t$ et en utilisant le théorème de changement de base pour le foncteur ${Rf{\sR}}_!$, on obtient le diagramme commutatif 
\[\begin{tikzcd}
H^i_{c, \sC}(\sX_{\sP}, \sF) & H^i_{c, \sC}(\sX_{\sR}, \sF)  \\
H^i_{c, \sP}(\sX_{\sP}, \sF) & H^i_{c, \sR}(\sX_{\sR}, \sF)  \\
	\arrow[from=1-1, to=1-2]
	\arrow[from=1-2, to=2-2]
	\arrow[from=1-1, to=2-1]
    \arrow[from=2-2, to=2-1]
\end{tikzcd}\]

D'après le lemme \ref{lemnulQ(j)}, le groupe $H^{2d+3-a}_{c, \sC}(\sX_{\sP}, \hat{\sT} \otimes \Q(d+1-a))$ est nul. On en déduit le diagramme

\[\begin{tikzcd}[row sep={40,between origins}, column sep={80,between origins}]
      & H^{2d+2-a}_{c, \sC}(\sX_{\sP}, \tilde{\sT}_t) \ar{rr}\ar{dd}\ar{dl} & & H^{2d+2-a}_{c, \sC}(\sX_{\sR}, \tilde{\sT}_t)\vphantom{\times_{S_1}} \ar{dd}\ar{dl} \\
    H^{2d+2-a}_{c, \sP}(\sX_{\sP}, \tilde{\sT}_t) \ar{dd} & & H^{2d+2-a}_{c, \sR}(\sX_{\sR}, \tilde{\sT}_t) \ar[crossing over]{ll} \\
      & H^{2d+2-a}_{c, \sC}(\sX_{\sP}, \tilde{\sT})  \ar{rr} \ar{dl} \ar{dd}& &  H^{2d+2-a}_{c, \sC}(\sX_{\sR}, \tilde{\sT})\vphantom{\times_{S_1}} \ar{dl} \\
     H^{2d+2-a}_{c, \sP}(\sX_{\sP}, \tilde{\sT}) && H^{2d+2-a}_{c, \sR}(\sX_{\sR}, \tilde{\sT}) \ar[from=uu,crossing over] \ar[crossing over]{ll} \\
     & 0. \\
\end{tikzcd}\]

Puisque le carré horizontal supérieur et les carrés verticaux sont commutatifs, une chasse au diagramme permet d'obtenir la commutativité voulue. \qed \\

Les groupes $\mathcal{D}^{2d+2-a}_{c}(\sX_\sR,\tilde{\sT})$ vont nous permettre d'approximer dans un certain sens le groupe $\Sha^{2d+2-a}_c(\sX, \sS, \tilde{\sT})$; c'est l'objet de la proposition suivante.
\begin{prop}\label{propstabilisationD}
    Pour tout nombre premier $l \neq p$, il existe un ouvert non vide $\sR_0$ de $\sU$ contenant $\sS$ tel que pour tout ouvert $\sR \subseteq \sR_0$ contenant $\sS$
    $$ \mathcal{D}^{2d+2-a}_c(\sX_\sR,\tilde{\sT})\{l\}=\mathcal{D}^{2d+2-a}_c(\sX_{\sR_0},\tilde{\sT})\{l\}=\Sha^{2d+2-a}_{c}(\sX,\sS,\tilde{\sT})\{l\}$$
\end{prop}

\proof
Le lemme 3.7 de \cite{HarSza16} indique que toute suite décroissante de groupes abéliens de torsion $l$-primaire et de type cofini stabilise.
Or, pour tout ouvert $\sP \subseteq \sR$, on a que $ \mathcal{D}^{2d+2-a}_c(\sX_{\sP},\tilde{\sT})\{l\} \subseteq \mathcal{D}^{2d+2-a}_c(\sX_{\sR},\tilde{\sT})\{l\}$ par la proposition \ref{propcomutativitédiagrammeD} et ces groupes sont de type cofini. Ainsi, il existe un ouvert $\sR_0$ de $\sU$ contenant $\sS$ tel que pour tout ouvert $\sR \subseteq \sR_0$ contenant $\sS$,

$$ \mathcal{D}^{2d+2-a}_c(\sX_\sR,\tilde{\sT})\{l\}=\mathcal{D}^{2d+2-a}_c(\sX_{\sR_0},\tilde{\sT})\{l\}.$$ 

Soit $\alpha \in \mathcal{D}^{2d+2-a}_c(\sX_{\sR_0},\tilde{\sT})\{l\}$ et soit $v \in S$. Considérons $\sR$ un ouvert de $\sR_0$ contenant $\sS$ mais ne contenant pas $v$ (i.e. $v \in R$), on peut alors voir $\alpha$ comme un élément de $\mathcal{D}^{2d+2-a}_c(\sX_{\sR},\tilde{\sT})\{l\}$ par ce qui précède. Alors, puisque le morphisme $H^{2d+2-a}_{c, \sR}(\sX_{\sR}, \tilde{\sT}) \ra H^{2d+2-a}_{c, k_v}(X_v, \tilde{\sT})$ se factorise par $H^{2d+2-a}_{c, \sR}(\sX_{\sR}, \tilde{\sT}) \ra H^{2d+2-a}_{c, \sS}(\sX,\tilde{\sT})$,  le lemme \ref{lemsuitelocalisation} montre que l'image de $\alpha $ dans $H^{2d+2-a}_{c, k_v}(X_v, \tilde{\sT})$ est nulle, d'où une première inclusion
$$
\mathcal{D}^{2d+2-a}(\sX_{\sR_0},\tilde{\sT})\{l\} \subseteq \Sha^{2d+2-a}_c(\sX, \sS,\tilde{\sT})\{l\}.
$$

Montrons maintenant l'autre inclusion : soit $\beta \in \Sha^{2d+2-a}_c(\sX, \sS,\tilde{\sT}))$. On peut relever $\beta$ en un élément de $H^{2d+2-a}(\sX_{\sR}, \tilde{T})$ pour un certain $\sR$ contenant $\sS$ et comme l'image de $\beta$ dans $\bigoplus_{v \in R} H^{2d+2-a}_{c, k_v}(X_v, \tilde{\sT}))$ est nulle par définition de $\beta$, le lemme \ref{lemsuitelocalisation} montre que $\beta$ est l'image d'un élément de $H^{2d+2-a}_{c, \sC}(\sX_\sR,\tilde{\sT})$. On en déduit la seconde inclusion 
$$
\Sha^{2d+2-a}_c(\sX, \sS,\tilde{\sT})\{l\} \subseteq \mathcal{D}^{2d+2-a}(\sX_{\sR_0},\tilde{\sT})\{l\}.
$$

et donc l'égalité
$$
\Sha^{2d+2-a}_c(\sX, \sS,\tilde{\sT})\{l\} =\mathcal{D}^{2d+2-a}(\sX_{\sR_0},\tilde{\sT})\{l\}.
$$
\qed


\begin{lemme}\label{lemsuitelocal2}
 Soit $\sR$ un ouvert de $\sU$ contenant $\sS$ et soit $\sF$ un faisceau localement constant sur $\sR$. La suite suivante
\begin{equation}\label{suitelemHarari}
\bigoplus_{v \in S} H^{i-1}(k_v^h, \sF) \ra H^i_c(\sR, \sF) \ra H^i(\sS, \sF)
\end{equation}
est exacte.
\end{lemme}

\proof La preuve est strictement identique à celle de la proposition 4.2  de \cite{HarSza16} après avoir remplacé le point générique par $\sS$, on en donne ici les grandes lignes. 
Soit $\alpha \in \bigoplus_{v \in S} H^{i-1}(k_v^h, \sF)$ et soit $\sP$ un ouvert de $\sR$ contenant $\sS$ choisi de sorte que $\alpha \in \bigoplus_{v \in P} H^{i-1}(k_v^h, \sF).$

On considère le diagramme commutatif suivant à ligne exacte 
\[\begin{tikzcd}
	\bigoplus_{v \in P} H^i(k_v,\sF) & H^{i}_{c}(\sP, \sF) & H^{i}(\sP, \sF) \\
	&  H^{i}_{c}(\sR,\sF) & H^{i}(\sR, \sF).
	\arrow[from=1-1, to=1-2]
	\arrow[from=2-2, to=2-3]
	\arrow[from=1-2, to=1-3]
    \arrow[from=1-2, to=2-2]
    \arrow[from=2-3, to=1-3]
\end{tikzcd}\]
Rappelons que le morphisme $H^i_c(\sR, \sF) \ra H^i(\sS, \sF)$ est  défini comme la composition $H^i_c(\sR, \sF) \ra H^{i}(\sR, \sF) \ra H^i(\sS, \sF)$. La suite horizontale du haut étant exacte, l'image de $\alpha$ dans $ H^{i}(\sP, \sF)$ est nulle; il en va donc de même de son image dans $H^i(\sS, \sF)$. Par commutativité du diagramme, l'image de $\alpha$ dans $H^{i}_{c}(\sR,\sF)$ est également envoyée sur l'élément nul de $H^i(\sS, \sF)$; cela prouve que la suite \ref{suitelemHarari} est un complexe.\\

Pour montrer que ce complexe est bien exact, considérons le diagramme commutatif à ligne exacte (\cite{HarSza16}, proposition 3.1,2) :
\[\begin{tikzcd}
	 H^i_c(\sP,\sF) & H^{i}_{c}(\sR, \sF) & \displaystyle \bigoplus_{v \in P \setminus R}H^{i}(\kappa(v), \sF) \\
	&  H^{i}(\sS,\sF) & \displaystyle \bigoplus_{v \in P \setminus R}H^{i}(k_v, \sF)
	\arrow[from=1-1, to=1-2]
	\arrow[from=2-2, to=2-3]
	\arrow[from=1-2, to=1-3]
    \arrow[from=1-2, to=2-2]
    \arrow[from=1-3, to=2-3]
\end{tikzcd}\]
où le morphisme vertical de droite est défini pour toute place $v \in \sR$ par $H^{i}(\kappa(v), \sF) \simeq H^{i}(\CO_v^h, \sF) \ra H^{i}(k_v^h, \sF) \ra H^{i}(k_v, \sF)$. Soit maintenant $\beta \in Ker\bigg(H^{i}_{c}(\sR, \sF) \ra H^{i}(\sS,\sF) \bigg)$. À fortiori, l'image de $\beta$ dans $\displaystyle \bigoplus_{v \in \sR \setminus \sP}H^{i}(k_v, \sF)$ est également nulle. De plus, le morphisme $H^i(\mathcal{O}_v^h, \sF) \ra H^{i}(K_v^h, \sF)$ est injectif d'après le lemme 1.3 de \cite{Sa89} ce qui implique que le morphisme $H^{i}(\kappa(v), \sF) \ra H^{i}(k_v, \sF)$ l'est également. Par conséquent, l'image de $\beta$ dans $\bigoplus_{v \in \sR \setminus \sP}H^{i}(\kappa(v), \sF)$ est nulle et étant donné que la suite horizontale du haut est exacte (\cite{HarSza16}, proposition 3.1, (2)), on peut relever $\beta$ en un élément de $H^i_c(\sP,\sF)$. Puisque $\displaystyle H^{i}(\sS,\sF) \coloneqq \drl_{\sP \subseteq \sU, \sU \setminus \sP \subseteq S} H^i(\sP, \sF)$, on peut quitte à rétrécir $\sP$ supposer que l'image de $\beta$ dans $H^i(\sP, \sF)$ est triviale. L'exactitude de la suite $\bigoplus_{v \in P} H^i(k_v,\sF) \ra H^{i}_{c}(\sP, \sF) \ra H^{i}(\sP, \sF)$ permet alors de conclure. \qed \\

On adapte ensuite ce résultat au cadre qui nous intéresse.

\begin{prop}\label{propsuiteexacteDc}
 Soit $\sR $ un ouvert de $\sU$ contenant $\sS$. La suite
$$
\bigoplus_{v \in S}H^{2d+1-a}_{c,k_v}(X_v ,\tilde{\sT}) \ra H^{2d+2-a}_{c, \sC}(\sX_{\sR} , \tilde{\sT}) \ra \mathcal{D}^{2d+2-a}_{c}(\sX_{\sR},\tilde{\sT}) \ra 0
$$
est exacte.
    
\end{prop}

\proof

Le lemme \ref{lemsuitelocal2} fournit pour tout faisceau $\sF$ localement constant sur $\sR$ une suite exacte 
$$ \bigoplus_{v \in S} H^{i-1}(k_v^h, \sF) \ra H^i_c(\sR, \sF) \ra H^i(\sS, \sF).$$\\

Prenons $\sF = {Rf_{\sR}}_{!}(\hat{\sT} \otimes \Z / n\Z(d+1-a))$ avec $n$ premier à $p$, ce faisceau est bien localement constant.
En appliquant le théorème de changement de base, on obtient une suite exacte

$$
 \bigoplus_{v \in S} H^{2d+1-a}_{c, k_v}(\sX \otimes_{k} k_v^h,\hat{\sT} \otimes \Z / n\Z(d+1-a)) \ra H^{2d+2-a}_{c, \sC}(\sX_\sR,\hat{\sT} \otimes \Z / n\Z(d+1-a)) \ra \mathcal{D}^{2d+2-a}_c(\sX_{\sR},\hat{\sT} \otimes \Z / n\Z(d+1-a)) \ra 0$$
qui en passant à la limite inductive devient
$$ 
\bigoplus_{v \in S} H^{2d+1-a}_{c, k_v}(\sX \otimes_{k} k_v^h,\hat{\sT} \otimes \Q / \Z(d+1-a)) \ra H^{2d+2-a}_{c, k_v}(\sX_\sR,\hat{\sT} \otimes \Q / \Z(d+1-a)) \ra \mathcal{D}^{2d+2-a}_c(\sX_\sR,\hat{\sT} \otimes \Q / \Z(d+1-a)) \ra 0. $$ \\

D'après la proposition 10.1 de \cite{Ge17}, le morphisme naturel $H^{2d+1-a}_{c, k_v^h}(\sX \otimes_{k} k_v^h,\hat{\sT} \otimes \Q / \Z(d+1-a)) \ra H^{2d+1-a}_{c, k_v}(X_v,\hat{\sT} \otimes \Q / \Z(d+1-a))$ est un isomorphisme. Le lemme \ref{lemnulQ(j)} donne la nullité des groupes \begin{itemize}
    \item $H^{2d+3-a}_{c, \sC}(\sX_{\sR}, \hat{\sT} \otimes \Q(d+1-a)),$
    \item $H^{2d+3-a}_{c, \sR}(\sX_{\sR}, \hat{\sT} \otimes \Q(d+1-a)),$
    \item $H^{2d+2-a}_{c, k_v}(X_v,  \hat{\sT} \otimes \Q(d+1-a)),$
    \item $H^{2d+3-a}_{c, k_v}(X_v, \hat{\sT} \otimes \Q(d+1-a)),$
\end{itemize}
ce qui nous permet de considérer le diagramme commutatif suivant dont la première ligne et les colonnes sont exactes

\[\begin{tikzcd}
	&& 0 \\
	{ \bigoplus_{v \in S} H^{2d+1-a}_{c,k_v}(X_v,\tilde{\sT}_t)} & {H^{2d+2-a}_{c,\sC}(\sX_\sR,\tilde{\sT}_t)} & {\mathcal{D}^{2d+2-a}_{c}(\sX_\sR,\tilde{\sT}_t)} & {0} \\
	{ \bigoplus_{v \in S} H^{2d+1-a}_{c, k_v}(X_v,\tilde{\sT})} & {H^{2d+2-a}_{c, \sC}(\sX_\sR,\tilde{\sT})} & {\mathcal{D}^{2d+2-a}_c(\sX_\sR,\tilde{\sT})} & {0} \\
	0 & 0 & 0.
	\arrow[from=2-2, to=2-3]
    \arrow[from=3-1, to=4-1]
	\arrow[from=2-1, to=2-2]
	\arrow[from=2-1, to=3-1]
	\arrow[from=2-2, to=3-2]
	\arrow[from=3-2, to=4-2]
	\arrow[from=2-3, to=3-3]
    \arrow[from=2-3, to=2-4]
    \arrow[from=3-3, to=3-4]
	\arrow[from=3-2, to=3-3]
	\arrow[from=3-1, to=3-2]
	\arrow[from=3-3, to=4-3]
	\arrow[from=1-3, to=2-3]
\end{tikzcd}\]

Une chasse au diagramme permet alors de déduire l'exactitude de la seconde ligne de la première. \qed \\

Nous aurons besoin du petit lemme suivant.

\begin{lemme}\label{lemdivisible}
   Soit $A$ un groupe abélien et $l$ un nombre premier. Alors, $\ol{A}\{l\} \simeq \ol{A\{l\}}$. 
\end{lemme}

\proof Comme les groupes divisibles sont injectifs dans la catégorie des groupes abéliens, on a en notant $Div_{max}$ le sous-groupe divisible maximal de $A$ une suite exacte scindée
$$
0 \ra Div_{max} \ra A \ra \ol{A} \ra 0.
$$
Ainsi, $A \simeq Div_{max} \bigoplus \ol{A}$ et on obtient le résutat voulu en prenant la partie $l$-primaire puis en quotientant par le sous-groupe divisible maximal de chaque côté.\qed \\

\begin{thme}\label{thmedualitégénérale}
Supposons que le groupe $H^{a+1}(X_v,\Z (a) )\{l\}$ est d'exposant fini pour tout $v \in S$ et tout nombre premier $l$. 
On a un accouplement parfait de groupes finis
    $$
    \ol{\Sha^{a}(\sX,\sS,\sT)}\{p'\} \times \ol{\Sha^{2d+2-a}_c(\sX,\sS,\tilde{\sT})}\{p'\} \ra \Q / \Z\{p'\} 
    $$
\end{thme}

\proof

 Pour tout ouvert $\sR \subseteq \sU$ contenant $\sS$,on définit le groupe $D^a(\sX_\sR, \sT)$ via la suite exacte 

$$ 0 \ra D^a(\sX_\sR, \sT) \ra H^{a}( \sX_\sR, \sT) \ra \prod_{v \in S} H^a( X_v , \sT).$$

Pour toute place $v \in S$, le groupe $H^{a+1}( X_v , \Z(a))\{l\}$ est d'exposant fini par hypothèse; par restriction-corestriction, il en est de même du groupe $H^a( X_v , \sT)\{l\}$. Ce groupe n'admet donc pas de sous-groupe divisible non trivial. En outre, étant donné que le groupe $H^a(\sX_\sR,\check{\sT} \otimes \Q(a))$ est sans torsion, la suite
$$
H^a(\sX_\sR,\check{\sT} \otimes \Q / \Z(a)) \ra  H^a(\sX_\sR,\sT) \ra H^a(\sX_\sR,\check{\sT} \otimes \Q(a))
$$
combiné au lemme \ref{lemtypecofini} montre que le groupe $H^a(\sX_\sR,\sT)\{l\}$ est de torsion et de type cofini. Il en résulte par passage au quotient des sous-groupes divisibles maximaux la suite exacte suivante 

$$ 0 \ra \overline{D^a(\sX_\sR, \sT)\{l\}} \ra \overline{H^a(\sX_\sR,\sT)\{l\}} \ra \prod_{v \in S} H^a( X_v , \sT)\{l\}. $$

De plus, on a d'après la proposition \ref{propstabilisationD} la suite exacte 
$$
\bigoplus_{v \in S}H^{2d+1-a}_{c, k_v}(\sX\otimes_{\sS} k_v,\tilde{\sT}) \ra H^{2d+2-a}_{c, \sC}(\sX_\sR , \tilde{\sT}) \ra \mathcal{D}^{2d+2-a}_c(\sX_\sR ,\tilde{\sT}) \ra 0
$$
 et comme les groupes $H^{2d+2-a}_{c, \sC}(\sX_\sR , \tilde{\sT})\{l\}$ et $\bigoplus_{v \in S}H^{2d+1-a}_{c, k_v}(X_v,\tilde{\sT})\{l\}$ sont de torsion de type cofini d'après le lemme \ref{lemtypecofini}, il en découle une suite exacte 

 $$
 \bigoplus_{v \in S}\overline{H^{2d+1-a}_{c, k_v}(\sX\otimes_{\sS} k_v,\tilde{\sT})}\{l\} \ra \overline{H^{2d+2-a}_{c, \sC}(\sX_\sR , \tilde{\sT})}\{l\} \ra \overline{\mathcal{D}^{2d+2-a}_c(\sX_\sR ,\tilde{\sT})}\{l\}
 \ra 0
 $$
 en vertu du lemme \ref{lemdivisible}.\\
Ces deux suites forment le diagramme commutatif à lignes exactes suivant

\[\begin{tikzcd}
	0 & {\overline{D^a(\sX_\sR, \sT)\{l\}}} & {\ol{H^a(\sX_\sR,\sT)\{l\}}} & {\prod_{v \in S} H^{a}(X_v,\sT)\{l\}} \\
	0 & (\overline{\mathcal{D}^{2d+2-a}_c(\sX_\sR ,\tilde{\sT})}\{l\})^D& (\overline{H^{2d+2-a}_{c, \sC}(\sX_\sR , \tilde{\sT})}\{l\})^D & (\bigoplus_{v \in S}\overline{H^{2d+1-a}_{c, k_v}(\sX\otimes_{\sS} k_v,\tilde{\sT})}\{l\})^D .
	\arrow[from=1-2, to=1-3]
	\arrow[from=1-1, to=1-2]
	\arrow[from=2-1, to=2-2]
	\arrow[from=2-2, to=2-3]
	\arrow[from=2-3, to=2-4]
	\arrow[from=1-2, to=2-2]
	\arrow[from=1-4, to=2-4]
	\arrow[from=1-3, to=2-3]
	\arrow[from=1-3, to=1-4]
\end{tikzcd}\]

Puisque $H^a(\sX_\sR,\sT)\{l\}$ est de type cofini, on a un isomorphisme $\ol{H^a(\sX_\sR,\sT)\{l\}} \simeq (H^a(\sX_\sR,\sT)\{l\})^{(l)}$. Les morphismes verticaux du milieu et de droite sont des isomorphismes par les propositions \ref{propdualitéglobale} et \ref{propdualitélocale}; une chasse au diagramme donne alors un isomorphisme $ D^{a}(\sX_\sR ,\sT)\{l\} \simeq (\overline{\mathcal{D}^{2d+2-a}_c(\sX_\sR ,\tilde{\sT})}\{l\})^D$. D'une part, on a pour un ouvert $\sR_0 \subseteq \sU$ suffisamment petit que $\mathcal{D}_c^{2d+2-a}(\sX_{\sR_0} ,\tilde{\sT})\{l\} = \Sha_c^{2d+2-a}(\sX,\sS,\tilde{\sT})\{l\}$ grâce à la proposition \ref{propstabilisationD}.
En passant à la limite directe sur les $\sR \subseteq \sR_0$ contenant $\sS$, on peut conclure que $$\ol{\Sha^{a}(\sX,\sS,\sT)}\{p'\} \simeq (\overline{\Sha^{2d+2-a}_c(\sX,\sS,\tilde{\sT})}\{p'\})^D$$.\\


\qed

\begin{rmke}
 Pour $l \neq p$, le groupe $H^{a+1}(X_v,\Z (a) )\{l\}$ est d'exposant fini lorsque $X_v$ est projectif sur $k_v$, a bonne réduction et que $a \neq 0, 1 $ ou $2$, cf (\cite{Ge17}, théorème 4.5).
\end{rmke}

\subsection*{Le cas des tores}
Rappelons que $\sT \coloneqq \check{\sT} \otimes \Z(a)[1]$. Lorsque $a=1$, on a que $\Z(1)[1] \simeq \gm$ et le faisceau $\sT$ est localement (pour la topologie étale) isomorphe à un produit fini de $\gm$, c'est donc un tore. Dans ce cadre, on peut affiner le théorème \ref{thmedualitégénérale} et établir la dualité suivante.\\

\begin{thme}
On a un accouplement parfait de groupes finis
    $$\Sha^{1}(\sX,\sS,\sT)\{p'\} \times \overline{\Sha^{2d+1}(\sX,\sS,\tilde{\sT})}\{p'\} \ra \Q / \Z\{p'\}.$$
\end{thme}

\proof Soit $v \in S$. Le groupe $Pic(X_v)$ s'inscrit dans une suite exacte courte
\begin{equation}\label{suiteexactePictores}
    0 \ra Pic^0(X_v) \ra Pic(X_v) \ra NS(X_v) \ra 0.
\end{equation} 

D'une part, le groupe $NS(X_v)$ est de type fini par le théorème de Néron-Severi (\cite{Né52}, théorème 2, page 145). D'autre part, le groupe $Pic^0(X_v)$ s'identifie aux $k_v$-points d'une variété abélienne définie sur le corps local $k_v$. On aura besoin du petit lemme suivant.

\begin{lemme}
Soit $k$ un corps local complet de caractéristique $p>0$ et notons $R$ l'anneau des entiers de $k$ et $k_0$ son corps résiduel. Soit $A$ une variété abélienne définie sur $k$. Alors, la torsion du groupe $A(k)$ est finie.
\end{lemme}

\proof
Soit $\sA$ le modèle de Néron de $k$ sur $R$. D'après la proposition 10 de \cite{ClXA08}, la torsion du noyau du morphisme de réduction $\phi : \sA(R) \ra \sA(k_0)$; on peut alors conclure grâce à la suite exacte
$$
0 \ra Ker(\phi) \ra \sA(R) \ra \sA(k_0) \ra 0
$$
puisque $\sA(R) \simeq A(k).$
\qed

Ainsi, la suite exacte (\ref{suiteexactePictores}) montre que la torsion de $\Pic(X_v)$ est finie. Par restriction-corestriction, on en déduit la finitude du groupe $H^2(X_v, \Z(1) )_{tors} \simeq H^1(X_v, \G_m )_{tors}$.\\



En outre, la proposition 6.1 de \cite{GuJaWi96} montre que le groupe de Picard de $\sX$ est de type fini. Cela implique que le groupe $H^1(\sX, \sT)\{l\}$ est d'exposant fini et n'a donc pas de sous-groupe divisible non-trivial. L'accouplement du théorème \ref{thmedualitégénérale} devient donc $$\Sha^{1}(\sX,\sS,\sT)\{p'\} \times \overline{\Sha^{2d+1}(\sX,\sS,\tilde{\sT})}\{p'\} \ra \Q / \Z\{p'\}.$$ \qed

\textbf{Remerciements :} Je remercie Diego Izquierdo pour son importante aide dans la réalisation et dans l'écriture de ce article;
son soutien a été essentiel afin de mener ce travail à bien.

\bibliographystyle{alpha}
\bibliography{Biblio}
\end{document}